\documentclass{amsart}[11pt]
\usepackage{amssymb}
\usepackage{verbatim}
\usepackage{pb-diagram}
\usepackage{hyperref}

\usepackage[margin=1.35in]{geometry}

\newtheorem{exercise}{Exercise}
\newtheorem{lem}[exercise]{Lemma}
\newtheorem{prop}[exercise]{Proposition}

\newtheorem{theorem}[exercise]{Theorem}
\newtheorem{conj}[exercise]{Conjecture}

\theoremstyle{definition}

\newtheorem{remark}[exercise]{Remark}

\newcommand{\A}{\mathbb{A}}
\newcommand{\AAA}{\mathcal{A}}
\newcommand{\Ad}{\mathrm{Ad}}
\newcommand{\ad}{\mathrm{ad}}
\newcommand{\af}{\mathrm{af}}
\newcommand{\al}{\alpha}

\newcommand{\B}{\mathbb{B}}

\def\C{{\mathbb C}}

\newcommand{\Des}{\mathrm{Des}}
\newcommand{\ee}{\hat{e}}

\newcommand{\F}{\mathbb{F}}
\newcommand{\fin}{\mathrm{fin}}

\newcommand{\Frac}{\mathrm{Frac}}
\newcommand{\geh}{{\mathfrak{g}}}
\newcommand{\gr}{j}
\newcommand{\Gr}{\mathrm{Gr}}
\newcommand{\he}{\hat e}
\newcommand{\hh}{\hat h}
\newcommand{\hhh}{\mathfrak{h}}

\newcommand{\hL}{\hat\Lambda}
\newcommand{\hs}{\hat s}

\newcommand{\id}{\mathrm{id}}
\newcommand{\ip}[2]{\langle #1\,,\,#2 \rangle}
\newcommand{\la}{\lambda}
\newcommand{\La}{\Lambda}

\newcommand{\pt}{\mathrm{pt}}
\newcommand{\Pet}{{\mathbb P}}
\newcommand{\Q}{\mathbb{Q}}
\newcommand{\qSchub}{\tilde{\mathfrak{S}}}

\def\Schub{{\mathfrak{S}}}

\newcommand{\Sym}{\mathrm{Sym}}
\newcommand{\tphi}{\tilde \phi}
\newcommand{\tpsi}{\tilde \psi}

\newcommand{\tS}{\tilde{S}}
\newcommand{\We}{W_e}
\newcommand{\Y}{\mathbb{Y}}
\def\Z{{\mathbb Z}}

\title[Quantum double Schuberts to $k$-double Schurs]
{From double quantum Schubert polynomials to $k$-double Schur functions via the Toda lattice}

\author{Thomas Lam}
\address{Department of Mathematics,
University of Michigan, 530 Church St., Ann Arbor, MI 48109 USA}
\email{tfylam@umich.edu}
 \thanks{T.L. was supported by NSF grant DMS-0901111, and by a Sloan Fellowship.}
\author{Mark Shimozono}
\address{Department of Mathematics, Virginia Tech, Blacksburg, VA 24061-0123 USA}
\email{mshimo@vt.edu}
\thanks{M.S. was supported by NSF DMS-0652641 and DMS-0652648.}
\begin{document}
\begin{abstract}
We show that the $k$-double Schur functions defined by the authors,
and the quantum double Schubert polynomials studied by Kirillov and Maeno
and by Ciocan-Fontanine and Fulton,
can be obtained from each other by an explicit rational substitution.  The main new ingredient is an explicit computation of Kostant's
solution to the Toda lattice in terms of equivariant Schubert
classes.
\end{abstract}

\maketitle


\section{Introduction}
Our purpose is to establish that the quantum double Schubert polynomials studied in \cite{KM} \cite{CF}
are related to the $k$-double Schur functions  of the authors \cite{LamSh:double} by an
explicit rational substitution of variables. 
The quantum double Schubert polynomials are the Schubert basis of the torus-equivariant quantum cohomology ring 
$QH^T(SL_n/B)$ of the flag manifold $SL_n/B$ \cite{AC} \cite{LamSh:qSchub}. The authors have shown \cite{LamSh:double} that the $k$-double Schur functions 
are a symmetric function realization (see the map $\kappa$ in \eqref{E:thediagram})
of the equivariant Schubert basis of the equivariant homology ring $H_T(\Gr_{SL_n})$ where $\Gr_{SL_n}$ 
is the affine Grassmannian of $SL_n$.

The context for our work is the following commutative diagram of isomorphisms, which is explained in the
remainder of this introductory section. The isomorphisms going around the square are valid
for any simply-connected semisimple algebraic group $G$.

\begin{align}\label{E:thediagram}
\begin{diagram}
\node{QH^T(G/B)_{(q)}} \arrow{e,t}{\psi}  \arrow{s,t}{\rho} \node{H_T(\Gr_G)_{\mathcal{T}}} \arrow{e,t}{\kappa} 
\node{\La_{(n)}(y\|a)_\bullet} \\
\node{\C[\mathcal{A}^\circ\times_{\hhh/W} \hhh]} \arrow{e,b}{\Psi} \node{\C[Z^\circ]} \arrow{n,b}{\phi} 
\end{diagram}
\end{align}
The vertical isomorphisms may be viewed as providing explicit presentations
for the top rings as coordinate rings of varieties. The map $\rho$ is Kim's presentation of $QH^T(G/B)$ \cite{Kim}
and the map $\phi$ is an isomorphism due independently to Ginzburg \cite{Gin} and Peterson \cite{P}.  Peterson's quantum to affine isomorphism $\psi$ \cite{P,LamSh:Qaff,LL} is defined in terms of Schubert classes,
while the map $\Psi$, which is Kostant's solution to the Toda lattice \cite{Ko:toda}, is defined in terms of explicit rational substitutions of coordinates.

Forgetting equivariance, we recover the setting of \cite{LamSh:qSchub}.

Our proof consists of computing images of variables in Kim's presentation
under Kostant's isomorphism, and writing these elements in terms of affine Schubert classes.  By Kostant's work \cite{Ko:toda,Kos}, the images of these generators are ratios of minors of the matrix entries of a centralizer family $Z$. 
In the nonequivariant setting \cite{LamSh:qSchub} these minors are directly recognizable as Schubert classes
in $\Gr_{SL_n}$, but equivariantly this does not hold. However the minors are
Schubert classes in the larger ring $H_T(\Gr_{PGL_n})$ and as such, have a nice
determinantal formula (Theorem \ref{T:adjointcentHom} and Proposition \ref{P:mapdet}).

Peterson \cite{P} has given a model for $H_T(\Gr_G)$, where $G$ is simply-connected and semisimple, in terms of Kostant and Kumar's nilHecke ring \cite{KK}.
To compute in $H_T(\Gr_{G_\ad})$ where $G_\ad$ is the adjoint group of the same type as the simply-connected group $G$, we use Chaput, Manivel, and Perrin's generalization \cite{CMP} of Peterson's construction to an extended affine nilHecke construction of $H_T(\Gr_{G_\ad})$.
Computation of Schubert structure constants  in the extended affine setting, appear to yield new positivity properties (Conjecture \ref{con:Mih}), despite the lack of a geometric interpretation coming from the positivity of the product in $QH^T(G/B)$ \cite{Mih}.

Finally,  we compare the minors of the centralizer with symmetric functions, using Molev's Jacobi-Trudi formula \cite{Mo} for dual Schur functions, and show that the quantum double Schubert polynomials correspond to $k$-double Schur functions
up to a twist by an automorphism of the symmetric function ring.

\subsection{Peterson's quantum versus affine Theorem}
We review Peterson's results relating quantum and affine Schubert calculus \cite{P}, following the presentation in \cite{LamSh:Qaff}.  See also \cite{LL}.

Let $G\supset B\supset T$ be a semisimple simply-connected algebraic group $G$ over $\C$,
a Borel subgroup $B$, and maximal torus $T$, $W=N(T)/T$ the Weyl group, and $I$ the set of Dynkin nodes.

Let $QH^T(G/B)$ be the $T$-equivariant small quantum cohomology ring of $G/B$ \cite{Kim}.
It contains a polynomial subring generated by $H_2(G/B)\cong Q^\vee$ where $Q^\vee$
is the coroot lattice. $H_2(G/B)$ has basis given by the
quantum parameters $\{q_i\mid i\in I\}$; $q_i$ corresponds to the simple coroot $\alpha_i^\vee$.
Let $S=H^T(\pt)$. $QH^T(G/B)$ has a basis over $S[q_i\mid i\in I]$
given by the quantum equivariant Schubert classes $\sigma^w$ for $w\in W$. 

On the other hand, let $\Gr_G=G(\C((t)))/G(\C[[t]]))$ be the affine Grassmannian of $G$.
Since $\Gr_G$ is an affine Kac-Moody homogeneous space, its equivariant homology 
$H_T(\Gr_G)$ is a free $S$-module with equivariant Schubert basis
$\xi_w$ for $w\in W_\af^0$ where $W_\af$ is the affine Weyl group and $W_\af^0$ is the
set of minimum length coset representatives in $W_\af/W$ \cite{KK}. $H_T(\Gr_G)$ also carries an
$S$-linear Pontryagin product since $\Gr_G$ is weakly homotopy equivalent to the space of based loops
into the compact form of $G$ and the latter space has a $T$-equivariant Pontraygin product
coming from the product in the target of loops.

We have $W_\af \cong W \ltimes Q^\vee$ where $\mu\mapsto t_\mu$ denotes the 
embedding $Q^\vee\to W_\af$. For $\la\in Q^\vee$ antidominant, $t_\la\in W_\af^0$.
The Schubert basis has the factorization property 
\begin{align}\label{E:factQ}
  \xi_{w t_\mu} = \xi_w \xi_{t_\mu}\qquad\text{for $w\in W_\af^0$ and $\mu\in Q^\vee$ antidominant.}
\end{align}
By \eqref{E:factQ} the set
\begin{align}\label{E:transclasses}
\mathcal{T}=\{\xi_{t_\la}\mid \text{$\la\in Q^\vee$ is antidominant}\}\subset H_T(\Gr)
\end{align}
is multiplicatively closed. 
Finally, for any $w\in W_\af$ there is a sufficiently antidominant $\la\in Q^\vee$
such that $wt_\la\in W_\af^0$.

\begin{theorem} \label{T:quantum2affine} \cite{P} \cite{LamSh:Qaff,LL} 
There is an $S$-algebra isomorphism
\begin{align*}
  \psi: QH^T(G/B)_q &\cong H_T(\Gr_G)_{\mathcal{T}} \\
  \sigma^w q_\mu &\mapsto \xi_{w t_\la} \xi_{t_{\la-\mu}}^{-1}
\end{align*}
for all $w\in W$ and $\mu\in Q^\vee$ and any $\la\in Q^\vee$ antidominant
such that $wt_\la\in W_\af^0$ and $\la-\mu$ is antidominant,
where $QH^T(G/B)$ is localized at the quantum parameters and $H_T(\Gr_G)$ is localized
at the set $\mathcal{T}$.
\end{theorem}

This is the map $\psi$ in \eqref{E:thediagram}.

\subsection{Kim's presentation} 
\label{SS:GK}

We recall a presentation \cite{GK} \cite{Kim} for $QH^T(SL_n/B)$ or more generally $QH^T(G/B)$.
Consider the $n\times n$ matrix $C_n$ with diagonal entries $x_i$ for $1\le i\le n$,
superdiagonal entries all $-1$, and subdiagonal entries $q_i$ for $1\le i\le n-1$.
\begin{align}\label{E:tridi}
C_n=  \begin{pmatrix}
  x_1 & -1 & \\
  q_1 & x_2 & -1 \\
      & q_2 & x_3 & -1 \\
&&\ddots&\ddots&\ddots \\
&&&q_{n-2}&x_{n-1}& -1 \\
 && & & q_{n-1}&x_n
  \end{pmatrix}
\end{align}
Let $g_{j,n}(x;q)\in\C[x;q]$ be the value of the $j$-th basic invariant evaluated at $C_n$:
\begin{align}\label{E:intmotion}
  \det(C_n-z\,\id_n) &= \sum_{j=0}^n (-z)^{n-j} g_{j,n}(x;q).
\end{align}
We let $S \cong \C[a_1,\dotsc,a_n]/(a_1+\dotsm+a_n)$; the $a_i$ are coordinates
on the Cartan subalgebra of $\mathfrak{sl_n}$.

Let $J$ be the ideal in $S[x;q]$ generated by $g_{j,n}(x;q)-e_j(a)$
for $1\le j\le n$ where $e_j(a)=e_j(a_1,\dotsc,a_n)$ is the elementary symmetric polynomial.

\begin{theorem} \label{T:QHT} \cite{GK} \cite{Kim} There is an $S$-algebra isomorphism
\begin{align*}
  \rho=\rho_{SL_n}:QH^T(SL_n/B) \cong S[x;q]/J.
\end{align*}
\end{theorem}

For any simply-connected semisimple algebraic group $G$ 
we recall a formulation of Kim's construction of
$QH^T(G/B)$ (see also \cite{R}).
Let $\geh=\mathrm{Lie}(G)$ with Cartan subalgebra $\hhh=\mathrm{Lie}(T)\subset \geh$.
Let $G^\vee$ be the Langlands dual group, $\geh^\vee$ its Lie algebra,
and $\hhh^\vee$ the Cartan subalgebra. The complexified root lattice of $\geh$ (resp. $\geh^\vee$)
is isomorphic to $\hhh^\vee$ (resp. $\hhh$). A perfect pairing of the root lattices of $\geh^\vee$
and $\geh$ is given by $\ip{\alpha_i^\vee}{\alpha_j}=a_{ij}$ where $(a_{ij}\mid i,j\in I)$ is the Cartan matrix.
By duality this induces a perfect pairing $\hhh \times \hhh^\vee\to \C$. There is an embedding $\hhh\to\geh^{\vee*}$
such that $h\in\hhh$ maps to the functional on $\geh^\vee$ that vanishes on all the root subspaces of $\geh^\vee$,
and on $\hhh^\vee$ is given by $h'\mapsto \ip{h}{h'}$ for $h'\in\hhh^\vee$.

$G^\vee$ acts on $\geh^{\vee*}$ by the coadjoint action
\begin{align*}
  (g\cdot f)(x) = f(\Ad(g^{-1})\cdot x)\qquad\text{for $g\in G^\vee$, $f\in \geh^{\vee*}$ and $x\in \geh^\vee$.}
\end{align*}
Let $E = \sum_{i=1}^{n-1} f_i^{\vee*}\in \geh^{\vee*}$ where $f_i^{\vee*}\in \geh^{\vee*}$ is the functional
taking the value $1$ on the basis element $f_i^\vee$ of the one-dimensional root subspace $\geh^\vee_{-\alpha^\vee_i}$
and zero on other weight spaces. Let $e_i^{\vee*}\in \geh^{\vee*}$ be defined similarly. 

Consider the affine scheme and open subscheme
\begin{align}
\AAA &= (-E+\hhh) \oplus \bigoplus_{i\in I} \C e_i^{\vee*}\subset \geh^{\vee*} \\
\AAA^\circ &= (-E+\hhh) \oplus \bigoplus_{i\in I} \C^\times e_i^{\vee*}\subset \geh^{\vee*}
\end{align}
There is a morphism $\AAA\to \hhh/W$ where $W$ is the Weyl group,
defined by the inclusion $\AAA \subset \geh^{\vee*}$ followed by the quotient map
$\geh^{\vee*} \to \geh^{\vee*}/G^\vee\cong \hhh/W$. 
Then Kim's theorem \cite{Kim} is that there is an $S$-algebra isomorphism
\begin{align}
  \rho: QH^T(G/B) \cong \C[\AAA \times_{\hhh/W} \hhh].
\end{align}

Let us consider the above construction for $\geh=\mathfrak{sl}_n(\C)$. We have $G^\vee=PGL_n$.
Consider the pairing $\mathfrak{gl}_n \times \mathfrak{gl}_n\to \C$ given by
$(A,B)\mapsto \mathrm{trace}(AB)$. This induces a perfect pairing
$\mathfrak{sl}_n \times \mathfrak{pgl}_n \to \C$ which extends
$\hhh\times\hhh^\vee\to \C$. The pairing gives an isomorphism $\mathfrak{pgl}_n^* \cong \mathfrak{sl}_n$,
under which $e_i^{\vee*}$ (resp. $f_i^{\vee*}$) is identified with the matrix in $\mathfrak{sl}_n$ having entry $1$
in position $(i+1,i)$ (resp. $(i,i+1)$) and zeroes elsewhere. 
The presentation of Theorem \ref{T:QHT} is recovered by translating
the general construction for this special case, where the $x$ variables are coordinates
on $\hhh$ in $\AAA$, the $q$ variables are the coordinates of the $e_i^{\vee*}$,
and the $a$ variables are the coordinates of the base scheme $\hhh$,
the right hand copy in $\AAA \times_{\hhh/W} \hhh$.

Localizing the isomorphism of Theorem \ref{T:QHT} at the quantum parameters
yields the map $\rho$ in \eqref{E:thediagram} where $\AAA^\circ\subset\AAA$ is the Zariski open
subset on which the $q_i$ are nonvanishing.

\subsection{Quantum double Schubert polynomials}
For $i\in I$ let $s_i^a$ be the operator that exchanges $a_i$ and $a_{i+1}$
and $\partial_i^a = (a_i-a_{i+1})^{-1} (1-s_i^a)$ the divided difference operator.
For $w\in S_n$ the quantum double Schubert polynomial
$\qSchub_w(x;q)\in S[x;q]$ is defined by \cite{KM, CF} \footnote{Our quantum double Schubert polynomials 
$\qSchub_wf(x;a)$ differ from those in \cite{KM} by negating the variables $a_i$.}
\begin{align*}
\qSchub_{w_0}(x;a) &= \prod_{i=1}^{n-1} \det(C_i-a_{n-i}\id_i) \\
\qSchub_w(x;a) &= -\partial_i^a \qSchub_{s_iw }(x;a)\qquad\text{if $s_iw>w$.}
\end{align*}

\begin{theorem} \label{T:QHTSchub} \cite{AC} \cite{LamSh:qSchub}
Under the isomorphism $\rho_{SL_n}$ of Theorem \ref{T:QHT}, 
for every $w\in S_n$ the $w$-th equivariant quantum Schubert class $\sigma^w\in QH^T(SL_n/B)$,
maps to the coset of the quantum double Schubert polynomial $\qSchub_w(x;a)$.
\end{theorem}

\subsection{Centralizer family of Ginzburg and Peterson}
\label{SS:centralizer}
For a semisimple algebraic group $G$,
following Peterson \cite{P} (cf. \cite{R}) let
\begin{align*}
  Z &= \{(h,g)\in \hhh\times G^\vee\mid g \cdot (-E+h) = -E+h\} \\
  &= \{(h,b)\in \hhh\times B^\vee \mid b\cdot (-E+h)=-E+h\}
\end{align*}
where $B^\vee\subset G^\vee$ is the Borel in the Langlands dual group $G^\vee$.

\begin{theorem} \label{T:centralizer} \cite{Gin} \cite{P}
There is an $S$-Hopf algebra isomorphism
\begin{align}
  \phi: \C[Z] \to H_T(\Gr_G).
\end{align}
\end{theorem}

Let $Z^\circ\subset Z$ be the open subset given by the
non-vanishing locus of the set $\mathcal{T}$ in \eqref{E:transclasses}. Then
\begin{align} \label{E:centralizerHom}
  \C[Z^\circ]\cong H_T(\Gr_G)_{\mathcal{T}}.
\end{align}
Alternatively, $(h,b) \in Z^\circ$ if and only if $b \in B^\vee_- w_0 B^\vee_-$ is in the open opposite Bruhat cell.  The localized isomorphism furnishes the map $\phi$ of \eqref{E:thediagram} when $G$ is simply-connected.

In \S \ref{S:SLcentralizer} we shall give an explicit description (and a new proof) of Theorem \ref{T:centralizer} in terms of Schubert classes for both the simply-connected group $G=SL_n$ and 
the adjoint group $G=PGL_n$.

\subsection{Symmetric function realization of $H_T(\Gr_{SL_n})$}
In this section $G=SL_n$. For historical reasons we use the notation $k=n-1$ here.
In \cite{LamSh:double} the $S$-Hopf algebra $H_T(\Gr_{SL_n})$ and its
equivariant Schubert basis is realized by symmetric functions.

\begin{theorem} \label{T:GrSym} \cite{LamSh:double} 
There is an $S$-Hopf algebra $\La_{(n)}(y\|a)$ of symmetric functions
and an $S$-basis consisting of the $k$-double Schur functions $s_\la^{(k)}(y\|a)$ and an isomorphism
of $S$-Hopf algebras
\begin{align}
  H_T(\Gr_{SL_n}) &\to \La_{(n)}(y\|a)
\end{align}
sending the Schubert basis to the $k$-double Schur basis. 
\end{theorem}

This theorem is explained, and made more precise, in \S \ref{SS:kdouble}.
After localizing this yields the map $\kappa$ in \eqref{E:thediagram}.  The $k$-double Schur functions are double analogues \cite{Mo} of Lapointe, Lascoux, and Morse's $k$-Schur functions \cite{LLM}.

\subsection{Kostant's isomorphism}
\label{SS:Kostant}

In Kostant's solution of the generalized Toda lattice he showed:

\begin{theorem} \label{T:Kostant} \cite[Thm. 2.4]{Ko:toda} 
There is an isomorphism
$Z^\circ\to \AAA^\circ\times_{\hhh/W} \hhh$ over $\hhh$. 
\end{theorem}
Denote by $\Psi$ (see \eqref{E:thediagram}) the induced isomorphism
\begin{align}\label{E:Psi}
  \Psi: \C[\AAA^\circ \times_{\hhh/W} \hhh] \to\C[Z^\circ].
\end{align}

Kostant also computes the images of various elements under this map.
We recall this as well as the explicit description for $Z^\circ$
in the case $G=SL_n$, which we assume for the rest of the section.
We have $G^\vee=PGL_n$. Let $\{z_{ij}\mid 1\le i\le j\le n\}$ be the projective matrix entry 
coordinates on $B^\vee$ and let $S=\C[\hhh]=\C[a_1,\dotsc,a_n]/(a_1+\dotsm+a_n)$.
Using the various identifications at the end of \S \ref{SS:GK}, we see that
$(h,b)\in Z$ if and only if 
\begin{align}\label{E:commeqs}
  z_{ij} = z_{i-1,j-1} + (a_{i-1}-a_j) z_{i-1,j}\qquad\text{for $2\le i\le j\le n$.}
\end{align}

Fix $k$ with $0\le k\le n$.
For a partition $\la=(\la_1,\la_2,\dotsc,\la_k)$ with $\la_1 \le n-k$ let
\begin{align}\label{E:zdet}
  z_{\la,k} = \det (z_{p,\la_{k+1-q}+q})_{1\le p,q\le k}.
\end{align}
It is the minor of the matrix $(z_{ij})$ using the first $k$ rows and the
columns given by the set of indices $(\la_k+1,\la_{k-1}+2,\dotsc,\la_1+k)$.

Let $R_k$ be the rectangular partition with $n-k$ rows and $k$ columns.
For any $0\le i\le k$ let $R_k-i$ be the partition $R_k$ but with $i$ cells
removed from the last column. 

Let $D_i=z_{R_i,n-i}$ and $D_i'=z_{R_i-1,n-i}$ for $0\le i\le n$.
Let $Z^\circ\subset Z$ be the Zariski open subset 
on which all the $D_i$ are nonvanishing.

Since $z_{ii}\ne0$ for all $1\le i\le n$ we may use the affine coordinates
$y_{ij}=z_{ij}/z_{11}$ for $1\le i\le j\le n$. In particular we let $y_{11}=1$.
Let $y_{\la,k}$, $D_i^y$, and $(D_i')^y$ be the quantities where $y$ variables are used instead of $z$ variables.
The following is the coordinatized version of Theorem \ref{T:Kostant} for $G=SL_n$.

\begin{theorem} \label{T:AKostant} \cite{Ko:toda} \cite{Kos} There is an
$S$-algebra isomorphism 
\begin{align}
 (S[x;q]/I)[q_1^{-1},\dotsc,q_{n-1}^{-1}]\cong   S[y_{ij}][y_{22}^{-1},\dotsc,y_{nn}^{-1},(D_0^y)^{-1},\dotsc,(D_{n-1}^y)^{-1}]
\end{align}
defined by
\begin{align}\label{E:trididiagonal}
  x_1+\dotsm+x_i &\mapsto a_1+\dotsm+a_i + \dfrac{D_i'}{D_i} \\
\label{E:tridisubdiagonal}
  q_i &\mapsto \dfrac{D_{i-1}D_{i+1}}{D_i^2}.
\end{align}
\end{theorem}

Note that $D_i'/D_i$ and $D_{i-1}D_{i+1}/D_i^2$ are equal to the analogous ratios involving $y$ variables;
in fact any ratio $f/g$ of homogeneous polynomials $f$ and $g$ in $z_{ij}$
of the same degree, is equal to the corresponding ratio in the $y$ variables.

\begin{remark} In \cite[Thm. 37, Cor. 28]{Kos} the 
respective formulae \eqref{E:trididiagonal} and \eqref{E:tridisubdiagonal}
are given for the zero fiber of the centralizer family, which is the nonequivariant case.
However these computations work for the general fiber
as well, although it requires additional work to compare these expressions with equivariant affine Schubert classes, as we shall do.
\end{remark}

\subsection{Main theorem}
Recall from \S \ref{SS:Kostant} the definition of the partitions $R_i$ and $R_i'=R_i-1$.
Let $(R_i)^t$ and $(R_i')^t$ be the transpose partitions. Let $\Des(w)=\{i\in I\mid ws_i<w\}$.

\begin{theorem} \label{T:main}
The rational transformation
\begin{align}
  x_1+\dotsm+x_i &\mapsto a_1+\dotsm+a_i + \dfrac{\hs_{(R'_i)^t}(y\|a)}{\hs_{(R_i)^t}(y\|a)} \\
  q_i &\mapsto \dfrac{\hs_{(R_{i-1})^t}(y\|a)\hs_{(R_{i+1})^t}(y\|a)}{\hs_{(R_i)^t}(y\|a)^2}.
\end{align}
sends
\begin{align}
\qSchub_w(x;a)\mapsto \frac{s^{(k)}_{\la(w)}(y\|a)^{\omega\eta}}{\prod_{i \in \Des(w)}\hs_{(R_i)^t}(y\|a)}
\end{align}
where $\la(w)$ is as described in \cite{LamSh:qSchub} and the automorphisms $\eta$ and $\tau$
are defined in \S \ref{SS:symfunc}. 
\end{theorem}

\begin{remark}
Forgetting equivariance, the map of Theorem \ref{T:main} is related to
the isomorphism of \cite{LamSh:qSchub} by the ``transpose'' automorphism $\omega$.
The map $\omega$ does not induce an isomorphism equivariantly. So perhaps a more geometrically
correct definition of a usual $k$-Schur function would incorporate either $\omega$
or the induced involution on the index set (induced by the affine Dynkin automorphism
$i\mapsto -i \mod n$). In fact one sees this ``transposed''
indexing in two geometric computations of homology Schubert classes in $\Gr_{SL_n}$ \cite{M} \cite{BS}.
\end{remark}

\subsection{Comments}
In principle, an analogue of Theorem \ref{T:main} holds for any semisimple group $G$, and one could try to use it to find symmetric function realizations of the Schubert basis of
$H_T(\Gr_G)$, from the corresponding quantum double Schubert polynomials, or vice versa.  However, it appears that neither the symmetric functions, nor the quantum double Schubert polynomials are known for any $G$ outside of type $A$.

In the nonequivariant setting, quantum Schubert polynomials were constructed in general type in \cite{Mar}.  Symmetric functions representing non-equivariant affine Schubert classes were constructed in \cite{Lam:Schub} for $SL_n$, \cite{LSS} for $Sp_{2n}$, and \cite{Pon}
for $SO_n$. See also Magyar's construction \cite{M}, which in principle can be applied for any $G$.  

In \cite{LamSh:qSchub} parabolic quantum double Schubert polynomials were defined and showed
to represent Schubert classes in $QH^T(SL_n/P)$ for a parabolic subgroup $P\subset SL_n$.
Moreover there is a parabolic analogue of Theorem \ref{T:quantum2affine}. One might try to
obtain a parabolic analogue of Theorem \ref{T:main}.

\subsection{Acknowledgements}
Thanks to Linda Chen for communicating to us her joint work with Dave Anderson \cite{AC}.
Thanks also to T. Ikeda for discussions about the Toda lattice.

\section{Extended affine nilHecke and Peterson algebras}
Let $G$ be simply-connected and let $G_\ad$ be the group of adjoint type with the same Lie algebra. 
Let $T\subset G$ be a maximal torus and let $T$ act on $G_\ad$ via $T\subset G\twoheadrightarrow G_\ad$.
In this section we recall the extended Peterson algebra, an algebraic model for $H_T(\Gr_{G_\ad})$
first studied by Chaput, Manivel, and Perrin \cite{CMP}.

\subsection{Extended affine Weyl group}
\label{SS:exaff}
Let $W_e = W \ltimes P^\vee$ be the extended affine Weyl group, where
$P^\vee=\bigoplus_{i\in I}\Z\omega_i^\vee$ is the coweight lattice and $\omega_i^\vee$ are the fundamental
coweights. Let $\Sigma'$ be the group of automorphisms of the affine
Dynkin diagram. Let $I_\af=I\cup\{0\}$ be the affine Dynkin node set.
Let $I^s = \Sigma' \cdot \{0\} \subset I_\af$ be the set of special nodes.
We have $P^\vee/Q^\vee = \{\omega_i^\vee+Q^\vee\mid i\in I^s\}$ where $\omega_0^\vee:=0$.
For every $i\in I^s$, the map $\omega_j^\vee+Q^\vee\mapsto \omega_j^\vee-\omega_i^\vee+Q^\vee$,
defines a permutation of the set $I^s$ which extends uniquely to an element
$\tau_i\in \Sigma'$ called the special automorphism associated with $i\in I^s$. It satisfies $\tau_i(i)=0$.
For $\mu\in P$ let $w_0^\mu\in W$ be the shortest element such that $w_0^\mu\cdot \mu = w_0\cdot \mu$.
We have
\begin{align*}
  \tau_i = w_0^{\omega_i} t_{-\omega_i^\vee}\qquad\text{for $i\in I^s$.}
\end{align*}
The special automorphisms $\Sigma =\{\tau_i\mid i\in I^s\}$ form a subgroup of $\Sigma'$.
We have $W_e = \Sigma \ltimes W_\af$
where $\sigma\in \Sigma$ acts on $W_\af$ by the group automorphism denoted 
$w\mapsto w^\sigma$ where $s_i^\sigma = s_{\sigma(i)}$ for all $i\in I_\af$. 
Moreover $\We$ admits a length function extending that of $W_\af$ and
$\Sigma$ is the subgroup of length-zero elements in $\We$.

$\We$ acts on the weight lattice $P$ of $G$, the polynomial ring $S=\Sym(P)$, and the fraction field $\F=\Frac(S)$ 
via the level zero action. Explicitly, the translation elements $t_\la$ for $\la\in P^\vee$ act on $P$ trivially.
$\sigma\in\Sigma$ acts on the affine root lattice $\bigoplus_{i\in I_\af} \Z \alpha_i$ by
$\alpha_i^\sigma = \alpha_{i^\sigma}$ for $i\in I_\af$ where
the superscript indicates the action of the automorphism $\sigma$. This induces an action
of $\Sigma$ on the quotient of the affine root lattice by the one-dimensional subspace
$\Z\delta$ where $\delta$ is the null root, and this quotient is isomorphic to the finite root lattice.
This action of $\Sigma$ extends uniquely to $P$.

\subsection{Extended small torus affine nilHecke algebra}
Let $\A_{\F} = \F\otimes_\Q \Q[\We]$ be the twisted group algebra. It has multiplication
$(p \otimes v)(q\otimes w) = p (v\cdot q)\otimes vw$ for $p,q\in \F$ and $v,w\in \We$.
The ring $\A_\F$ acts on $\F$: the field $\F$ acts on itself by left multiplication and $\We$ acts on $\F$ as above.

For $i\in I_\af$ define the divided difference element $A_i\in \A_\F$ by
\begin{align}\label{E:A}
  A_i = \alpha_i^{-1}(s_i-1).
\end{align}
Since the $A_i$ satisfy the same braid relations as $W_\af$,
for $v\in W_\af$ one may define $A_v = A_{s_{i_1}}\dotsm A_{s_{i_\ell}}$
where $v=s_{i_1}\dotsm s_{i_\ell}$ is a reduced decomposition.
For $\sigma\in \Sigma$ and $v\in W_\af$ define $A_{\sigma v} = \sigma A_v\in \A_\F$.
In particular $A_\sigma=\sigma$.
Since $A_i^2=0$ for $i\in I_\af$, one may show that for all $u,v\in \We$,
\begin{align}
A_u A_v = \begin{cases}
A_{uv} &\text{if $\ell(u)+\ell(v)=\ell(uv)$} \\
0 &\text{otherwise.}
\end{cases}
\end{align}

The small torus extended affine nilCoxeter algebra $\A_0$ is the subring of $\A_\F$ generated by
$\{\A_i\mid i\in I_\af\}$ and $\Sigma$. We have $$\A_0 = \bigoplus_{w\in \We} \Z A_w.$$
The ring $\A_0$ acts on $S$: for $i\in I_\af$, $\la\in P$, and $s,s'\in S$ we have
\begin{align}
  A_i \cdot \la &= -\ip{\alpha_i^\vee}{\la} \\
  A_i \cdot (s s') &= (A_i \cdot s) s' + (s_i\cdot s)(A_i\cdot s').
\end{align}

The small-torus finite (resp. affine, resp. extended affine) nilHecke algebra $\A_\fin$ (resp. $\A$, resp. $\A_e$) 
is by definition the subring of $\A_\F$ generated by $S$ and $\{A_i\mid i\in I\}$ (resp. $\{A_i\mid i\in I_\af\}$,
resp. $\A_0$). The action of $\A_\F$ on $\F$ restricts to an action of $\A_e$ on $S$.
We have
\begin{align*}
\A_e &\cong \bigoplus_{w\in W_e} S A_w.
\end{align*}
Note that there is an embedding $\We\to \A_e$ via
\begin{align}
  \sigma &\mapsto \sigma &\qquad&\text{for $\sigma\in \Sigma$} \\
  s_i &\mapsto 1+\alpha_i A_i &\qquad&\text{for $i\in I_\af$.}
\end{align}

\subsection{Coproduct}
For $S$-modules $M$ and $N$
define $M \otimes_S N$ to be the quotient of $M\otimes_\Z N$ by the submodule
generated by the elements $sm\otimes n - m\otimes sn$ for $s\in S$, $m\in M$, and $n\in N$.
Make a similar definition of $M\otimes_\F N$ for $\F$-modules $M$ and $N$.
Define the $\F$-submodule
\begin{align*}
  \Delta(\A_\F) = \bigoplus_{w\in \We} \F w \otimes w \subset \A_\F \otimes_\F \A_\F.
\end{align*}
The componentwise product on $\A_\F \otimes_{\Q} \A_\F$ induces an ill-defined product
on $\A_\F \otimes_{\F} \A_\F$ because elements of $\F$ do not
intertwine the same way with different elements of $\We$. However the componentwise product
does induce a well-defined product on $\Delta(\A_\F)$. This given, there is a ring 
and left $\F$-module homomorphism
$\Delta:\A_\F\to \Delta(\A_\F)$ given by
\begin{align}
  \Delta(q) &= q \otimes 1  &\qquad&\text{for $q\in \F$} \\
  \Delta(w) &= w \otimes w &&\text{for $w\in\We$.}
\end{align}
One may show that for all $i\in I_\af$,
\begin{align*}
  \Delta(A_i) &= A_i \otimes 1 + s_i \otimes A_i =  A_i \otimes 1 + 1 \otimes A_i + \alpha_i A_i \otimes A_i.
\end{align*}
Let $M$ and $N$ be $\A_\F$-modules. Then $M \otimes_\F N$ is an $\A_\F$-module
induced by the componentwise action of $\Delta(\A_\F)$ on $M \otimes_\F N$:
\begin{align}
  s \cdot (m \otimes n) &= \Delta(s) \cdot (m\otimes n) = (s \cdot m) \otimes n = m \otimes (s \cdot n) \\
  A_i \cdot (m \otimes n) &= \Delta(A_i) \cdot (m\otimes n) = (A_i \cdot m) \otimes n +( s_i \cdot m) \otimes (A_i \cdot n).
\end{align}
One may make analogous definitions for $\A_e$ instead of $\A_\F$.
For $\A_e$-modules $M$ and $N$, $\A_e$ acts on $M \otimes_S N$ as above.

\subsection{Extended Peterson subalgebra and $j$-basis}
The Peterson subalgebra $\B\subset\A$ and its extended analogue $\B_e\subset\A_e$ 
are by definition the centralizer subalgebras
\begin{align*}
\Pet &= Z_\A(S) = \{a\in \A \mid s a = a s\text{ for all $s\in S$ }\} \\
\Pet_e&=Z_{\A_e}(S)=\{a\in \A_e \mid s a = a s\text{ for all $s\in S$ }\}.
\end{align*}
We have
\begin{align*}
  \F \otimes_S \Pet &= \bigoplus_{\mu\in Q^\vee} \F \, t_\mu \\
  \F \otimes_S \Pet_e &= \bigoplus_{\mu\in P^\vee} \F \, t_\mu.
\end{align*}

Define the \textit{extended affine Grassmannian elements} in $\We$ to be the subset
\begin{align*}
\We^0 &= \{ w\in \We\mid \ell(wv)>\ell(w)\text{ for all $v\in W\setminus \id$}\}  \\
  &= \{\sigma w\mid \sigma\in \Sigma, w\in W_\af^0\}.
\end{align*}

\begin{lem} \label{L:Grassmannian} Let $w\in W$ and $\la\in P^\vee$. Then
$wt_\la\in \We^0$ if and only if $\ip{\la}{\alpha_i} \le0$ for all $i\in I$ and
strict inequality for $i\in I$ such that $ws_i<w$.
\end{lem}

The following result is the natural extension of a theorem of Peterson \cite{P} \cite{Lam:Schub} to the adjoint group setting, and its proof goes through the same way.  Parts (2) and (3) were established in \cite{CMP} and (1) is nearly immediate.  See also \cite{LSS} for a more algebraic approach.

\begin{theorem} \label{T:Petalg}
\begin{enumerate}
\item
The map $\Delta:\A_e\to\Delta(\A_e)$
restricts to an $S$-algebra homomorphism $\Pet_e \to \Pet_e\otimes_S \Pet_e$, making $\Pet_e$ into a
commutative $S$-Hopf algebra. 
\item
There is a unique $S$-basis $\{j_w\mid w\in W_e^0\}$ of $\Pet_e$ such that
\begin{align}
  j_w = A_w + \sum_{x\in W_e\setminus W_e^0} j_w^x A_x
\end{align}
for some elements $j_w^x\in S$. 
\item There is an $S$-Hopf algebra isomorphism 
\begin{align}
\label{E:toB}
\theta: H_T(\Gr_{G_\ad})&\cong \Pet_e \\
 \xi_w &\mapsto j_w
\end{align}
where $\{\xi_w\mid w\in \We^0\}$ is the Schubert basis of $H_T(\Gr_{G_\ad})$.
\end{enumerate}
\end{theorem}

Let $\gr:\A_e\to\Pet_e$ be the projection to the Peterson subalgebra; it is the left $S$-module
homomorphism defined by
\begin{align*}
  \gr(\sum_{w\in \We} a_w A_w) = \sum_{w\in \We^0} a_w j_w\qquad\text{for $a_w\in S$.}
\end{align*}
Alternatively, $j(a) \in \Pet_e$ is determined by the requirement that $j(a)-a \in  \sum_{v \in W \setminus \id} \A_eA_v$.

\begin{lem}\label{L:j} For all $a,a'\in \A_e$ and $b\in \Pet_e$,
\begin{align}
\label{E:jbb}
\gr(b) &= b \\
\label{E:jaa}
\gr(a'a) &= j(a'\gr(a)) \\
\label{E:jba}
\gr(ba) &= b \gr(a).
\end{align}
\end{lem}
\begin{proof} Equation \eqref{E:jbb} holds by the definitions.
$a-\gr(a)$ has no term $A_v$ for $v\in \We^0$; therefore the same is true of
$a'(a-\gr(a))$. Therefore $\gr(a'(a-\gr(a)))=0$ from which \eqref{E:jaa} follows.
Equation \eqref{E:jba} follows from \eqref{E:jbb} and \eqref{E:jaa}.
\end{proof}

\subsection{$j$-basis and Dynkin automorphisms}

\begin{lem}\label{L:conjnonextended} For any $\sigma\in\Sigma$,
conjugation by $\sigma$ defines a group automorphism of $W_\af$ and $W_e$
and a ring automorphism of $\A_e$, $\A$, $\Pet_e$, and $\Pet$.
\end{lem}

Let $i\mapsto i^*$ be the involutive automorphism of the affine Dynkin diagram
such that $w_0\cdot \alpha_i = -\alpha_{i^*}$ for $i\in I$ and $0^*=0$.
Then for $i\in I^s$ we have $\tau_i^{-1}=\tau_{i^*}$ and
$(w_0^{\omega_i})^{-1} = w_0^{\omega_{i^*}}$.

\begin{lem}\label{L:jaut}
\begin{align*}
  j_{\tau_i} = t_{\omega_{i^*}^\vee} \qquad\text{for $i\in I^s$.}
\end{align*}
\end{lem}
\begin{proof} 
We have
\begin{align*}
    t_{\omega_{i^*}^\vee} 
 = t_{-w_0\cdot\omega_i^\vee} 
 = w_0^{\omega_i} t_{-\omega_i^\vee} (w_0^{\omega_i})^{-1} 
 = \tau_i w_0^{\omega_{i^*}}.
\end{align*}
Since $w_0^{\omega_{i^*}}\in W$, its expansion in $\A_e$ into the $A$ basis,
has no Grassmannian terms other than $A_{\id}$. Therefore $t_{\omega_{i^*}^\vee}\in \Pet_e$
has no Grassmannian terms other than $A_{\tau_i}$. The lemma follows.
\end{proof}

For $a \in \Pet$ and $\tau \in \Sigma$, let $a^\tau:=\tau a \tau^{-1}$.

\begin{lem} \label{L:jtwist} For $w\in W_\af^0$ and $i\in I^s$ we have
\begin{align}\label{E:jtwist}
  j_{\tau_i w} = j_{\tau_i} j_w^{\tau_i}.
\end{align}
\end{lem}
\begin{proof} By the proof of Lemma \ref{L:jaut},
\begin{align*}
  j_w^{\tau_i} j_{\tau_i} &= j_w^{\tau_i} \tau_i w_0^{\omega_{i^*}^\vee} = \tau_i j_w w_0^{\omega_{i^*}^\vee}.
\end{align*}
Since $w_0^{\omega_{i^*}^\vee}\in W$, the right hand side is an element of $\Pet_e$
with unique Grassmannian term $A_{\tau_i w}$.
\end{proof}

\begin{lem}\label{L:jtranswt} Let $\la\in P^\vee$ be antidominant. Then
\begin{align}
  j_{t_\la} = \sum_{\mu\in W\cdot \la} A_{t_\mu}.
\end{align}
\end{lem}
\begin{proof} The proof is the same as that for the case that $\la\in Q^\vee$; see
\cite[Proposition 4.5]{Lam:Schub}. 
\end{proof}


\begin{lem} \label{L:Schubtrans}
Let $w\in \We^0$ and let $\mu\in P^\vee$ be antidominant. Then
$w t_\mu\in \We^0$ and
\begin{align}
  j_{wt_\mu} &= j_{t_\mu} j_w.
\end{align}
\end{lem}
\begin{proof} 
The only possible Grassmannian terms in the RHS are of the form $A_{t_{u\cdot \mu}} A_w$ where $u \in W$.  But this product is length additive only when it is equal to $A_{wt_\mu}$.
\end{proof}

\begin{remark}
Specializing to $G=SL_n$ and forgetting equivariance, the factorization of Lemma \ref{L:Schubtrans} is the $k$-rectangle factorization property for $k$-Schur functions \cite{LLM}.
\end{remark}

\subsection{Type $A$}
For type $A_{n-1}^{(1)}$, $I_\af=I^s=\Z/n\Z$, $(i+n\Z)^*= -i+n\Z$,
and $\tau_i=\tau^{-i}$ where $\tau$ is the rotation $j+n\Z\mapsto j+1+n\Z$.
$w_0^{\omega_i}\in W=S_n$
is the permutation sending $j\mapsto j+ n-i$ for $1\le j\le i$
and $j\mapsto j-i$ for $i+1\le j\le n$. Letting $u_i^j=s_i s_{i+1}\dotsm s_{i+j-1}$
and $d_i^j = s_{i+j-1} \dotsm s_{i+1} s_i = (u_i^j)^{-1}$ we have
we have $w_0^{\omega_i} = u_{n-i}^i \dotsm u_2^i u_1^i = d_1^{n-i} d_2^{n-i}\dotsm d_i^{n-i}$.
Lemma \ref{L:jaut} reads
\begin{align}\label{E:Ajrot}
  j_{\tau^k} &= t_{\omega_k^\vee}\qquad\text{for $0\le k\le n$}
\end{align}
with the convention that $\omega_0^\vee=\omega_n^\vee=0$.

\subsection{Positivity}
Define the subset of Graham-positive elements by $\Z_{\ge0}[\alpha_i\mid i\in I]\subset S$.
The coefficients $j^x_w$ are also structure constants of the $j_w$-basis. For the non-extended Peterson algebra they 
coincide with equivariant 3-point Gromov-Witten invariants, 
which are known to be Graham-positive by work of Mihalcea \cite{Mih}.

\begin{conj}\label{con:Mih}
All the coefficients $j^x_w$ are Graham-positive.
\end{conj}

\section{Symmetric functions}
The main purpose of this section is to establish Proposition \ref{P:detSchubclass}, which expresses certain ``small'' affine Schubert classes as a determinant of Dynkin-twists of special Schubert classes.

Let $\C[a_i\mid i\in\Z]$ be a polynomial ring over the complexes\footnote{It is possible to work over the integers.
However we use $\C$ since we work with coordinate rings of complex varieties.}, in variables indexed by integers.

\subsection{A Hopf algebra of symmetric functions}
\label{SS:symfunc}
Let $\hL(y\|a)$ be the $\C[a]$-Hopf algebra of symmetric series in the variables
$y=(y_1,y_2,\dotsc)$ with coefficients in $\C[a]$, with primitive elements
$p_r[y]=\sum_{i\ge1} y_i^r$ which generate $\hL(y\|a)$ in the sense that
$\hL(y\|a) = \C[a][[p_r[y]\mid r\ge1]]$.
$\hL(y\|a)$ is a completion of the usual $\C[a]$-algebra of symmetric functions
in the $y$ variables.

The dual elementary symmetric functions $\he_j(y\|a)\in\hL(y\|a)$ \cite{Mo} are defined by
\begin{align}\label{E:duale}
  \sum_{j\ge0} \he_j(y\|a) (t+a_0)(t+a_1)\dotsm(t+a_{j-1}) = \prod_i \dfrac{1+ty_i}{1-a_0y_i}.
\end{align}
Setting $t=-a_r$ for $r\in \Z_{\ge0}$ yields 
\begin{align}\label{E:dualerec}
  \sum_{j= 0}^r \he_j(y\|a) (a_0-a_r)(a_1-a_r)\dotsm(a_{j-1}-a_r) = \prod_i \dfrac{1-a_ry_i}{1-a_0y_i}.
\end{align}
One may solve for $\he_j(y\|a)$ in terms
of $\he_i(y\|a)$ for $i<j$; in particular the system \eqref{E:dualerec} uniquely defines the $\he_j(y\|a)$.
It also follows that $\he_j(y\|a)$ depends only on the parameters $a_0,a_1,\dotsc,a_j$.

Let $\tau$ be the $\C$-algebra automorphism of $\hL(y\|a)$ given by
\begin{align*}
a_i^\tau &= a_{i+1} &\qquad&\text{for $i\in\Z$} \\
p_r[y]^\tau &= p_r[y]&&\text{for $r\ge1$.}
\end{align*}

Let $\eta$ be the $\C$-algebra automorphism of $\C[a]$ defined by
\begin{align}
  a_i^\eta = - a_{1-i}\qquad\text{for $i\in\Z$.}
\end{align}
Let $\omega$ be the $\C[a]$-algebra automorphism of $\hL(y\|a)$ given by
\begin{align}
  s_\la[y]^\omega = s_{\la^t}[y]\qquad\text{for $\la\in\Y$}
\end{align}
where $\Y$ is Young's lattice of partitions and $\la\mapsto \la^t$
is the transpose or conjugate map.
Define the dual homogeneous symmetric functions $\hh_i(y\|a)\in\hL(y\|a)$ by
\begin{align}
  \hh_i(y\|a) = \ee_i(y\|a)^{\omega\eta}\qquad\text{for $i\ge0$.}
\end{align}
The dual Schur functions $\hs_\la(y\|a)\in\hL(y\|a)$ may be defined by \cite{Mo}
\begin{align}\label{E:dualSchur}
  \hs_\la(y\|a) &= \det (\hh^{\tau^{j-1}}_{\la_i-i+j}(y\|a))_{1\le i,j\le \ell(\la)} \\
  &= \det (\ee^{\tau^{1-j}}_{\la_i^t-i+j}(y\|a))_{1\le i,j\le \la_1}
\end{align}
Then we have
\begin{align}\label{E:etaversusomega}
  \hs_\la(y\|a)^{\omega\eta} = \hs_{\la^t}(y\|a) \qquad\text{for $\la\in\Y$.}
\end{align}


\subsection{A dual Hopf algebra of symmetric functions}
Let $\La(x\|a)$ be the $\C[a]$-Hopf algebra of symmetric functions
given by the polynomial ring over $\C[a]$ generated by the primitive elements
$p_r[x-a_+] = \sum_{i\ge1} (x_i^r-a_i^r)$ for $r\ge0$. Note that this coproduct
involves the $a_i$ variables nontrivially. We define a perfect pairing
$\ip{\cdot}{\cdot}: \La(x\|a)\times \hL(y\|a) \to \C[a]$ by
\begin{align}
\ip{p_\la[x-a_+]}{p_\mu[y]} = z_\la \delta_{\la,\mu}
\end{align}
where $z_\la = \prod_i i^{m_i} m_i!$ and $m_i$ is the number of times the part $i$ occurs in the partition $\la$.  We have \cite{Mo}
\begin{align*}
  \ip{s_\la(x\|a)}{\hs_\mu(y\|a)} = \delta_{\la\mu}
\end{align*}
where $s_\la(x\|a)$ is the double Schur function, which is essentially a limit involving
double Schubert polynomials indexed by Grassmannian permutations.

\begin{prop} \label{P:Hopfduality}
$\La(x\|a)$ and $\hL(y\|a)$ are Hopf dual over $\C[a]$ with respect to $\ip{\cdot}{\cdot}$.
In particular for $f,g\in\La(x\|a)$ and $h\in\hL(y\|a)$ we have
\begin{align}
\ip{f\otimes g}{\Delta(h)} &= \ip{fg}{h}.
\end{align}
\end{prop}

\subsection{The $k$-double Schur functions}
\label{SS:kdouble}
Let $k=n-1$. The ring $S=\C[a_1,\dotsc,a_n]/(a_1+a_2+\cdots+a_n)$ is a $\C[a]$-algebra
via the $\C$-algebra homomorphism $\C[a]\to S$ defined by
$a_{i+rn}\mapsto a_i$ for all $1\le i\le n$ and $r\in\Z$.

The $\C$-algebra automorphism $\eta$ of $\C[a]$ induces a
$\C$-algebra automorphism of $S$ (also denoted $\eta$) via the $\C[a]$-action on $S$.

In \cite{LamSh:double} the authors introduced a family of symmetric functions
$s_\la^{(k)}(y\|a)\in \hL(y\|a)$ for $\la\in\Y$ with $\la_1<n$ called $k$-double Schur functions.
They are linearly independent. 
Due to a nonstandard $S$-module structure used in \cite{LamSh:double},
we must apply the automorphism $\eta\circ\omega$ to revert to the standard $S$-action.

\begin{theorem} \label{T:kdoubleSchur} \cite{LamSh:double} 
Let $\La_{(n)}'(y\|a) = S \otimes_{\C[a]} \bigoplus_{\la_1<n} \C[a] s_\la^{(k)}(y\|a)^{\omega\eta}$.
Then 
\begin{enumerate}
\item
$\La_{(n)}'(y\|a)$ is an $S$-Hopf algebra with structure induced from $\hL(y\|a)$.
\item
The elements $1\otimes s_\la^{(k)}(y\|a)^{\omega\eta}$ are an $S$-basis of $\La_{(n)}'(y\|a)$.
\item
There is an $S$-Hopf isomorphism 
\begin{align*}
\kappa: H_T(\Gr_{SL_n})&\to \La_{(n)}'(y\|a) \\
\xi_{w_\la^\af}&\mapsto s_\la^{(k)}(y\|a)^{\omega\eta}
\end{align*}
where $\la\mapsto w_\la^\af$ is the bijection sending the $(n-1)$-bounded partition $\la\in\Y$ with $\la_1<n$, 
to $w_\la^\af\in W_\af^0$ (see \cite[Section 2.2]{LamSh:double}).  In particular, for $0\le r\le n-1$
we have $\xi_{c_r}\mapsto \hat{e}_r(y\|a)$ where $c_r=s_{r-1}\dotsm s_1s_0$.
\end{enumerate}
\end{theorem}

The effect of $S\otimes_{\C[a]} \cdot$
is to identify subscripts of the $a_i$ modulo $n$.  In \cite{LamSh:double}, we used the polynomial ring $\C[a_1,\ldots,a_n]$ instead of $S$, but the results there easily specialize ot the current situation.  The reader is also warned that the convention for simple roots in \cite{LamSh:double} 
is nonstandard: it uses $\alpha_i=a_{-i}-a_{1-i}$ for $i\in\Z$.
The automorphism $\eta$ is applied so that the standard $S$-actions can be used.
The automorphism $\omega$ is applied for convenience, 
so that small Schubert classes map to the known basis of dual Schur functions 
(but with transposed indexing partitions).

We denote by $\tau$ the special automorphism of the type $A_{n-1}^{(1)}$
affine Dynkin diagram, that sends $i$ to $i+1$ mod $n$ for all $i$. In the notation of
\S \ref{SS:exaff}, $\tau=\tau_1^{-1}$.

This automorphism acts on $H_T(\Gr_{SL_n})$ by conjugation.  In terms of the isomorphism $\Pet \simeq H_T(\Gr_{SL_n})$, $\tau$ acts via $a \mapsto a^\tau = \tau a \tau^{-1}$.

\begin{lem} \label{L:rotatehom} The isomorphism of Theorem \ref{T:kdoubleSchur}
intertwines the action of $\tau$ on $H_T(\Gr_{SL_n})$ and that induced by $\tau$
on $\La_{(n)}'(y\|a)$.
\end{lem}
\begin{proof}
In this proof we use notations from affine symmetric groups freely; see \cite[Section 2]{LamSh:double}.

Let $f \in \La(x\|a)$ and $g \in \La_{(n)}'(y\|a)$.  It is clear that $\ip{f}{\tau g} = \tau \ip{f}{g}$.  But there is also an evaluation map $\La(x\|a) \otimes \Pet \to S$, given by $\ip{f}{t_{\la}} = \epsilon_{t_\la}(f)$ (see \cite{LamSh:double}).  Here for $w \in W_\af = \tS_n$, the evaluation $\epsilon_w(f)$ is given by the substitution $x_i = -a_{w(1-i)}$ (this differs from the formula in \cite{LamSh:double} by $\eta$). It is enough to show that $\epsilon_{\tau t_\la \tau^{-1}}(f) = \tau \epsilon_{t_\la}(f)$ for any $\la \in Q^\vee$ and any $f \in \La(x\|a)$.  But $\tau t_{(\la_1,\la_2,\ldots,\la_n)} \tau^{-1} = t_{(\la_n,\la_1,\ldots,\la_{n-1})}$ so this follows from the proof of \cite[Lemma 9]{LamSh:double}, which states that $\epsilon_{t_{(\la_1,\la_2,\ldots,\la_n)}}(p_r[x-a_{>0}]) = -\sum_{i=1}^n \la_i a_i^r$. 
%
%
%
\end{proof}

\begin{lem} \label{L:smallkdoubleSchur}
For any $\la\in\Y$ with $\la_1\le n-k$ and $\ell(\la)\le k$ for some $1\le k\le n-1$ we have
\begin{align}
  s_\la^{(k)}(y\|a) = \hs_\la(y\|a).
\end{align}
\end{lem}
\begin{proof} 
We shall use the notations of \cite[Section 4]{LamSh:double}.  The partitions $\la$ of the lemma are exactly those with main hook length less than or equal to $n-1$.  Let $\mu$ have main hook length greater than $n-1$. 

By Lemma 9, the function $\epsilon_\Gr(s_\mu(x\|a))$ lies in the GKM ring $\Phi_\Gr$, and thus can be expanded in terms of $\xi^v$'s, where $v \in \tS_n^0$.  Suppose $\la$ has main hook length less than or equal to $n-1$.  Since $\la$ is not contained in $\mu$, we have $\epsilon_\la(s_\mu(x\|a)) = 0$ by \cite[Proposition 1]{LamSh:double}.
But since $\la$ is small, we also have $\epsilon_{w_\la^\af}(s_\mu(x\|a)) = 0$.  Thus the support of $\epsilon_\Gr(s_\mu(x\|a))$ does not contain any $w_\la^\af$ for $\la$ with main hook length less than or equal to $n-1$.

By duality, the expansion of $s^{(k)}_\la(y\|a)$ in terms of dual Schur functions does not involve $\hs_\mu(y\|a)$.  The lemma follows from this observation and \cite[Corollary 31]{LamSh:double}.
\end{proof}

\begin{prop} \label{P:detSchubclass}
For any $\la\in\Y$ with $\la_1\le n-k$ and $\ell(\la)\le k$ for some $1\le k\le n-1$ we have
\begin{align}
  \xi_{w_\la^\af} &= \det (\xi_{c_{\la_i-i+j}}^{\tau^{1-j}})_{1\le i,j\le k}.
\end{align}
\end{prop}
\begin{proof} Note that $\tau \eta = \eta \tau^{-1}$ and $\tau \omega = \omega\tau$. We have
\begin{align*}
  \kappa(\det(\xi_{c_{\la_i-i+j}}^{\tau^{1-j}})) 
  = \det(\ee_{\la_i-i+j}^{\tau^{1-j}}(y\|a)) 
  = \hs_{\la^t}(y\|a).
\end{align*}
On the other hand, we have
\begin{align*}
  \kappa(\xi_{w_\la^\af}) 
  = s_\la^{(k)}(y\|a)^{\omega\eta} 
  = \hs_\la(y\|a)^{\omega\eta} 
  = \hs_{\la^t}(y\|a)
\end{align*}
by Theorem \ref{T:kdoubleSchur}, Lemmata \ref{L:rotatehom} and \ref{L:smallkdoubleSchur}, and \eqref{E:etaversusomega}.
The Proposition follows since $\kappa$ is injective.
\end{proof}

\section{Centralizer family for $G^\vee=SL_n$}
\label{S:SLcentralizer}
Under the identification of $\mathfrak{pgl}_n^* \cong \mathfrak{sl}_n$ of \S \ref{SS:GK},
the element $E$ is mapped to the principal nilpotent (also denoted $E$) with entries $1$ on the superdiagonal
and zeroes elsewhere. With $\hhh$ still the Cartan subalgebra of $\mathfrak{sl}_n$, and using the adjoint action
of $SL_n$ on $\mathfrak{sl}_n$, define the family
\begin{align*}
 \tilde{Z} &= \{(h,b)\in \hhh\times SL_n \mid b\cdot (-E+h)=-E+h\} \\
 &= \{(h,b)\in \hhh \times B  \mid b\cdot (-E+h)=-E+h\}.
\end{align*}
The universal covering map $SL_n(\C) \to PGL_n(\C)$ induces a map $\tilde{Z}\to Z$ of schemes over $\hhh$
and an injective $S$-algebra homomorphism $\C[Z]\to \C[\tilde{Z}]$.

We denote the matrix entries of $B$ by $z_{ij}$ for $1\le i\le j\le n$.
The $z_{ij}$ satisfy $z_{ii}\ne0$ and $z_{11}z_{22}\dotsm z_{nn}=1$.  The ring $\C[\tilde{Z}]$ is generated by $S$ and $z_{ij}$, modulo the relations \eqref{E:commeqs}.

Let $c_p = s_{p-1} s_{p-2}\dotsm s_1 s_0\in W_\af$.

\begin{theorem} \label{T:adjointcentHom}
There is an $S$-Hopf algebra isomorphism $\tilde \phi:\C[\tilde{Z}]\to \Pet_e$ defined by
\begin{align}\label{E:phi}
  \tilde \phi(z_{k,k+p}) = \dfrac{j_{\tau^k c_p}}{j_{\tau^{k-1}}}
\end{align}
for $1\le k\le k+p\le n$.  Thus composing with $\theta^{-1}$ of Theorem \ref{T:Petalg}, we obtain the isomorphism $\phi: \C[\tilde Z] \to H_T(\Gr_{G_\ad})$.
\end{theorem}
\begin{proof} Let us first check well-definedness. Since the diagonal entries $z_{kk}$ are units in
$\C[\tilde{Z}]$ we must check that $\phi(z_{kk})$ is a unit in $\Pet_e$. Using \eqref{E:Ajrot} we have
\begin{align}\label{E:entries}
\tphi(z_{kk}) &= t_{\omega_k^\vee-\omega_{k-1}^\vee}\qquad\text{for $1\le k\le n$}.
\end{align}
But $t_\mu\in\Pet_e$ for all $\mu\in P^\vee$ and is invertible in $\Pet_e$.
The elements $j_{\tau^k}$ are invertible in $\Pet_e$ by 
\eqref{E:Ajrot}. This ensures that the division in the right hand side of \eqref{E:phi} is well-defined.

For $0\le k\le n$ we have
\begin{align}\label{E:phiomega}
  \tphi(z_{11}z_{22}\dotsm z_{kk}) = \prod_{i=1}^k t_{\omega_i^\vee-\omega_{i-1}^\vee} = t_{\omega_k^\vee}.
\end{align}
For $k=n$ this shows that the relation $z_{11}\dotsm z_{nn}-1$ is in the kernel of $\tphi$.

To check that the images of the equations \eqref{E:commeqs} are satisfied, proceeding by induction on $p$ 
it suffices to show
\begin{align}\label{E:goal}
  t_{-s_{k-2}\dotsm s_1 \omega_1^\vee} j_{\tau^k c_p} = j_{\tau^{k-1}c_p} + (a_{k-1}-a_{k+p}) j_{\tau^{k-1}c_{p+1}}.
\end{align}
We have
\begin{align*}
  t_{-\omega_1^\vee} &= (w_0^{\omega_1})^{-1} \tau_1 \\
  &= s_1 \dotsm s_{n-1} \tau^{-1} \\
  &= \tau^{-1} s_2\dotsm s_{n-1} s_0
\end{align*}
and
\begin{align*}
  t_{-s_{k-2}\dotsm s_1 \omega_1^\vee} &= s_{k-2}\dotsm s_1 t_{-\omega_1^\vee} s_1 \dotsm s_{k-2} \\
  &= s_{k-2}\dotsm s_1 (\tau^{-1} s_2\dotsm s_{n-1} s_0) s_1 \dotsm s_{k-2} \\
  &= \tau^{-1} s_{k-1}\dotsm s_2 s_2 \dotsm s_{n-1} s_0) s_1 \dotsm s_{k-2} \\
  &= \tau^{-1} s_k \dotsm s_{n-1} s_0 s_1 \dotsm s_{k-2}.
\end{align*}
Using Lemma \ref{L:j} repeatedly without further mention, we have
\begin{align*}
  t_{-s_{k-2}\dotsm s_1 \omega_1^\vee} j_{\tau^k c_p} 
  &= \gr(t_{-s_{k-2}\dotsm s_1 \omega_1^\vee} \tau^k A_{p-1}\dotsm A_1 A_0) \\
  &= \gr(\tau^{-1} s_k \dotsm s_{n-1} s_0 s_1\dotsm s_{k-2} \tau^k A_{p-1}\dotsm A_0) \\
  &= \gr(\tau^{k-1} s_0 \dotsm s_{n-2} A_{p-1}\dotsm A_0) \\
  &= \gr(\tau^{k-1} s_0 \dotsm s_{p-1} s_p A_{p-1}\dotsm A_0) \\
  &= \gr(\tau^{k-1} s_0\dotsm s_{p-1} (1+\alpha_p A_p) A_{p-1}\dotsm A_0) \\
  &= \gr(\tau^{k-1} s_0 \dotsm s_{p-1} A_{p-1}\dotsm A_0) + \gr(\tau^{k-1} s_0 \dotsm s_{p-1} \alpha_p A_p \dotsm A_0) \\
  &= \gr(\tau^{k-1} A_{p-1}\dotsm A_0) + \gr(\tau^{k-1} (\alpha_0+\dotsm+\alpha_p) s_0\dotsm s_{p-1} A_p \dotsm A_0) \\
  &= j_{\tau^{k-1}c_p} + (a_{k-1}-a_{k+p}) j_{\tau^{k-1}c_{p+1}}
\end{align*}
as required, since $\alpha_0+\dotsm+\alpha_p = a_0 - a_{p+1}$ and
applying $\tau^{k-1}$ yields $a_{k-1}-a_{p+k}$. 

Therefore $\tphi$ is well-defined. 

Next we check surjectivity.
To help with this we work with a localized map.
Let $\tilde{Z}^\circ$ be the locus in $\tilde{Z}$
on which the functions $D_i$ are nonvanishing. By Proposition \ref{P:mapdet}
$\tphi$ induces an $S$-algebra homomorphism
\begin{align} \label{E:philoc}
  \tphi_\bullet: \C[\tilde{Z}^\circ] \to (\Pet_e)_{(t)}
\end{align}
where $(\Pet_e)_{(t)}$ (resp. $\Pet_{(t)}$) is the localization of $\Pet_e$ (resp. $\Pet$) 
at the multiplicatively closed subset
consisting of the elements $j_{t_\la}$ for $\la\in P^\vee$ (resp. $\la\in Q^\vee$) antidominant.

Let $\tpsi:QH^T(G/B)_{(q)} \to \Pet_{(t)}$ be the composition of the map $\psi$ of Theorem \ref{T:quantum2affine}
with the isomorphism $H_T(\Gr_G)_{\mathcal{T}} \cong \Pet_t$ induced by the map \eqref{E:toB}.

For $i\in I$, by Lemma \ref{L:Grassmannian} we have $s_i t_{-\omega_i^\vee}\in \We^0$.
By the definition of $\psi$ and Lemma \ref{L:Schubtrans} we have
\begin{align*}
\tpsi(\sigma^{s_i}) &= j_{s_i t_{-\omega_i^\vee}} j_{t_{-\omega_i^\vee}}^{-1} \\
\tpsi(q_i) &= j_{t_{-\omega_{i-1}^\vee}} j_{-\omega_{i+1}^\vee} j^{-1}_{t_{-2\omega_i^\vee}}.
\end{align*}
By Proposition \ref{P:mapdet} we have
\begin{align*}
  \tphi(\Psi((x_1+\dotsm+x_i)-(a_1+\dotsm+a_i))) &=
  \tphi(D_i'D_i^{-1}) = \tpsi(\sigma^{s_i}) \\
  \tphi(\Psi(q_i)) &= \tphi(D_{i-1}D_{i+1}D_i^{-2}) = \tpsi(q_i).
\end{align*}
Since $\tpsi$ is an isomorphism and the elements
$q_i$ and $\sigma^{s_i}$ generate $QH^T(SL_n/B)_{(q)}$,
it follows that the image of $\tphi_{\bullet}$ contains $\Pet_{(t)}$.

By Lemmata \ref{L:jtwist} and \ref{L:jaut}, equation \eqref{E:phiomega}, and Lemma \ref{L:conjnonextended}
it follows that $\tphi_{\bullet}$ is surjective.
Using Lemma \ref{L:Schubtrans} we deduce that $\tphi$ is surjective.
Injectivity of $\tphi$ follows by dimension-counting.

To check that $\tphi$ is a coalgebra morphism, we note that $\Delta(z_{ij}) = \sum_{i\le k\le j} z_{ik}\otimes z_{kj}$.  In particular $\Delta(z_{ii}) = z_{ii} \otimes z_{ii}$ for $1\le i\le n$.  On the other hand, $\tphi(z_{ii})=t_{\omega_i^\vee-\omega_{i-1}^\vee}$ is grouplike.   Thus $\tphi$ is a coalgebra morphism.
\end{proof}

\begin{prop} \label{P:mapdet} 
\begin{align}\label{E:mapzdet}
  \tphi(z_{\la,k})= j_{\tau^k w_\la^\af}
\end{align}
for all partitions $\la=(\la_1,\dotsc,\la_k)$ with $\la_1 \le n-k$. In particular
\begin{align}
\label{E:Dtoj}
  \tphi(D_i) &= j_{t_{-\omega_i^\vee}} \\
\label{E:D'toj}
  \tphi(D'_i) &=  j_{s_i t_{-\omega_i^\vee}}
\end{align}
for $0\le i \le n-1$.
\end{prop}
\begin{proof} Applying Lemma \ref{L:jtwist}, reversing the rows and the columns of the
determinant and then transposing it, and using Proposition \ref{P:detSchubclass}  we have
\begin{align*}
  \tphi(z_{\la,k}) &= \det (j_{\tau^{p-1}}^{-1} j_{\tau^p c_{\la_{k+1-q}+q-p}}) \\
  &= (j_\tau j_{\tau^2} \dotsm j_{\tau^{p-1}})^{-1} \det (j_{\tau^p c_{\la_{k+1-q}+q-p}}) \\
  &= (j_\tau j_{\tau^2} \dotsm j_{\tau^{p-1}})^{-1} (j_\tau j_{\tau^2}\dotsm j_{\tau^k}) 
  \det (j_{c_{\la_{k+1-q}+q-p}}^{\tau^p}) \\
  &= j_{\tau^k} \det (j_{c_{\la_q-q+p}}^{\tau^{k+1-p}}) \\
  &= j_{\tau^k} (\det (j_{c_{\la_q-q+p}}^{\tau^{1-p}}))^{\tau^k} \\
  &= j_{\tau^k} j_{w_\la^\af}^{\tau^k} \\
  &= j_{\tau^k w_\la^\af}.
\end{align*}
Taking $\la=R_i$ with $n-i$ parts, we obtain 
\begin{align*}
  \tphi(D_i) &= j_{\tau^{n-i} w_{R_i}^\af} = j_{(w_{R_i}^\af)^{\tau^{n-i}} \tau^{n-i}} = j_{w_0^{\omega_{n-i}^\vee} \tau^{n-i}} = j_{t_{-\omega_i^\vee}}.
\end{align*}
For $D_i'$, by direct computation we have $w_{R_i'}^\af = s_{2i} w_{R_i}^\af$ with subscripts taken mod $n$. We have
\begin{align*}
  \tau^{n-i} w_{R_i'}^\af = \tau_{n-i} s_{2i} w_{R_i}^\af = s_i \tau^{n-i} w_{R_i}^\af = s_i t_{-\omega_i^\vee}
\end{align*}
as required.
\end{proof}

%

\section{Proof of Theorem \ref{T:main}}

\begin{proof}[Proof of Theorem \ref{T:main}]
The stated rational transformation is the composition of the maps $\Psi$, $\tilde \theta$,
the isomorphism $\theta^{-1}$ of Theorem \ref{T:Petalg}, and $\kappa$. 

The description of the map via the images of the $x$ and $q$ variables,
follows immediately from Theorems \ref{T:AKostant},
\ref{T:adjointcentHom}, \ref{T:Petalg}, \ref{T:kdoubleSchur}, Proposition \ref{P:mapdet},
and Lemma \ref{L:smallkdoubleSchur}.

The description in terms of Schubert classes follows immediately from Theorems
\ref{T:quantum2affine}, \ref{T:QHTSchub} and \ref{T:kdoubleSchur} together with the explicit correspondence between Schubert classes computed in \cite{LamSh:qSchub}.
\end{proof}

\end{document}